\documentclass[11pt]{article}
\usepackage{amssymb,amsmath,amsthm}
\newcommand{\R}{\mathbb{R}}
\newcommand{\Z}{\mathbb{Z}}
\newcommand{\N}{\mathbb{N}}
\def\qed{\hfill $\Box$ \smallskip}
\def\x#1{{\rm (\ref{#1})}}

\begin{document}
\title{\LARGE\bf{ Non-Nehari manifold method for asymptotically periodic Schr\"odinger equation}
 \footnote{This work is partially supported by the NNSF (No: 11171351) and the SRFDP(No: 20120162110021) of China;
 This paper has been accepted for publication in SCIENCE CHINA Mathematics. There are some errors in section 4 of the paper
 ``arXiv:1405.2607v1 [math.AP], 12 May 2014", this paper is also the revised version.}}
\date{}
 \author{ X. H. Tang\\
        {\small School of Mathematics and Statistics,}\\
        {\small Central South University,}\\
        {\small Changsha, Hunan 410083, P.R.China }\\
        {\small E-mail: tangxh@mail.csu.edu.cn}}
\maketitle
\begin{center}
\begin{minipage}{13cm}
\par
\small  {\bf Abstract:} We consider the semilinear Schr\"odinger equation
 $$
   \left\{
   \begin{array}{ll}
    -\triangle u+V(x)u=f(x, u), \ \ \ \  x\in {\R}^{N},\\
    u\in H^{1}({\R}^{N}),
   \end{array}
    \right.
 $$
 where $f$ is a superlinear, subcritical nonlinearity. We mainly study the case where
 $V(x)=V_0(x)+V_1(x)$, $V_0\in C(\R^N)$, $V_0(x)$ is 1-periodic in each of $x_1, x_2, \ldots, x_N$ and
 $\sup[\sigma(-\triangle +V_0)\cap (-\infty, 0)]<0<\inf[\sigma(-\triangle +V_0)\cap (0, \infty)]$,
 $V_1\in C(\R^N)$ and $\lim_{|x|\to\infty}V_1(x)=0$. Inspired by previous work of Li et al. \cite{LWZ},
 Pankov \cite{Pa} and Szulkin and Weth \cite{Sz}, we develop a more direct approach to
 generalize the main result in \cite{Sz} by removing the `` strictly increasing" condition in the
 Nehari type assumption on $f(x, t)/|t|$. Unlike the Nahari manifold method, the main idea of our
 approach lies on finding a minimizing Cerami sequence for the energy functional outside the
 Nehari-Pankov manifold $\mathcal{N}^{0}$ by using the diagonal method.


\end{minipage}
\end{center}

 {\section{Introduction}}
 \setcounter{equation}{0}

 \par
   Consider the following semilinear Schr\"odinger equation
 \begin{equation} \label{SP}
   \left\{
   \begin{array}{ll}
    -\triangle u+V(x)u=f(x, u), \ \ \ \  x\in {\R}^{N},\\
    u\in H^{1}(\R^N),
   \end{array}
    \right.
 \end{equation}
 where $V : {\R}^{N} \rightarrow {\R}$ and $f: {\R}^N\times {\R} \rightarrow {\R}$ satisfy the following
 basic assumptions, respectively:

 \vskip2mm
 \noindent
 \begin{itemize}
  \item[(V0)] $V\in C(\R^N)\cap L^{\infty}(\R^N)$, and
 \begin{equation}\label{sp}
  \sup[\sigma(-\triangle +V)\cap (-\infty, 0)]<0<\inf[\sigma(-\triangle +V)\cap (0, \infty)];
 \end{equation}
 \end{itemize}

 \par\noindent
 or

 \vskip2mm
 \noindent
 \begin{itemize}
  \item[(V0$'$)] $V(x)=V_0(x)+V_1(x)$, $V_0\in C(\R^N)\cap L^{\infty}(\R^N)$ and
 \begin{equation}\label{sp0}
  \sup[\sigma(-\triangle +V_0)\cap (-\infty, 0)]<0< \inf[\sigma(-\triangle +V_0)\cap (0, \infty)],
 \end{equation}
  $V_1\in C(\R^N)$ and $\lim_{|x|\to\infty}V_1(x)=0$;
 \end{itemize}

 \vskip2mm
 \noindent
 \begin{itemize}
  \item[(F0)] $f\in C(\R^N\times \R)$, and there exist constants $p\in (2, 2^*=2N/(N-2))$ and $C_0>0$ such that
 $$
   |f(x, t)|\le C_0\left(1+|t|^{p-1}\right), \ \ \ \ \forall \ (x, t)\in \R^N\times \R;
 $$

 \item[(F1)] $f(x, t)=o(|t|)$, as $|t|\to 0$, uniformly in $x\in \R^{N}$.
 \end{itemize}

 \vskip2mm
 \par
   Let $\mathcal{A}=-\triangle +V$ with $V\in C(\R^N)\cap L^{\infty}(\R^N)$. Then $\mathcal{A}$ is self-adjoint in
 $L^2(\R^N)$ with domain $\mathfrak{D}(\mathcal{A})=H^2(\R^N)$ (see \cite[Theorem 4.26]{EK}). Let $\{\mathcal{E}(\lambda):
 -\infty \le \lambda \le +\infty\}$ and $|\mathcal{A}|$ be the spectral family and the absolute value of $\mathcal{A}$,
 respectively, and $|\mathcal{A}|^{1/2}$ be the square root of $|\mathcal{A}|$. Set $U=id-\mathcal{E}(0)-\mathcal{E}(0-)$. Then $U$
 commutes with $\mathcal{A}$, $|\mathcal{A}|$ and $|\mathcal{A}|^{1/2}$, and $\mathcal{A} = U|\mathcal{A}|$ is
 the polar decomposition of $\mathcal{A}$ (see \cite[Theorem IV 3.3]{EE}). Let
 \begin{equation}\label{E}
   E=\mathfrak{D}(|\mathcal{A}|^{1/2}), \ \ \ \ E^{-}=\mathcal{E}(0-)E, \ \ \ \ E^{0}=[\mathcal{E}(0)-\mathcal{E}(0-)]E,
    \ \ \ \ E^{+}=[id-\mathcal{E}(0)]E.
 \end{equation}
 For any $u\in E$, it is easy to see that $u=u^{-}+u^{0}+u^{+}$ and
 \begin{equation}\label{Au}
   \mathcal{A}u^{0}=0, \ \ \ \ \mathcal{A}u^{-}=-|\mathcal{A}|u^{-}, \ \ \ \ \mathcal{A}u^{+}=|\mathcal{A}|u^{+},
       \ \ \ \ \forall \ u\in E\cap \mathfrak{D}(\mathcal{A}),
 \end{equation}
 where
 \begin{equation}\label{u-}
   u^{-}:=\mathcal{E}(0-)u\in E^{-}, \ \ \ \ u^{0}=[\mathcal{E}(0)-\mathcal{E}(0-)]u\in E^{0},
   \ \ \ \ u^{+}:=[id-\mathcal{E}(0)]u\in E^{+}.
 \end{equation}

 \par
    Note that $E^{0}=\mbox{Ker}(\mathcal{A})$, we can define an inner product
 \begin{equation} \label{IP}
    (u, v)=\left(|\mathcal{A}|^{1/2}u, |\mathcal{A}|^{1/2}v\right)_{L^2}+\left(u^{0}, v^{0}\right)_{L^2}, \ \ \ \ \forall \ u, v\in E
 \end{equation}
 and the corresponding norm
 \begin{equation} \label{No0}
    \|u\|=\left(\left\||\mathcal{A}|^{1/2}u\right\|_{2}^2+\left\|u^{0}\right\|_{2}^2\right)^{1/2}, \ \ \ \ \forall \ u\in E,
 \end{equation}
 where $(\cdot, \cdot)_{L^2}$ denotes the inner product of $L^2(\R^N)$, $\|\cdot\|_s$ denotes the norm of $L^s(\R^N)$.

 \par
    Firstly, we assume that (V0) is satisfied. In this case, $E^{0}=\{0\}$, $E=H^1(\R^N)$ with equivalent norms
 (see \cite{Di, DL}). Therefore, $E$ embeds continuously in $L^s(\R^N)$ for all $2\le s\le 2^*$. In addition, one
 has an orthogonal decomposition $ E=E^{-}\oplus E^{+}$, where orthogonality is with respect to both $(\cdot, \cdot)_{L^2}$
 and $(\cdot, \cdot)$. If $\sigma(-\triangle +V)\subset (0, \infty)$, then $E^{-}=\{0\}$, otherwise $E^{-}$ is infinite-dimensional.

 \par
   Under assumptions (V0), (F0) and (F1), the solutions of problem \x{SP} are critical points of the functional
 \begin{equation}\label{Ph}
   \Phi(u)=\frac{1}{2}\int_{{\R}^N}\left(|\nabla u|^2+V(x)u^2\right)\mathrm{d}x-\int_{{\R}^N}F(x, u)\mathrm{d}x,
 \end{equation}
 where $F(x, t)=\int_{0}^{t}f(x, s)\mathrm{d}s$. In view of \x{Au} and \x{No0}, we have
 \begin{equation}\label{Ph1}
   \Phi(u)  =   \frac{1}{2}\left(\|u^{+}\|^2-\|u^{-}\|^2\right)-\int_{{\R}^N}F(x, u)\mathrm{d}x,
                \ \ \ \ \forall \ u=u^{-}+u^{+}\in E^{-}\oplus E^{+}=E.
 \end{equation}

 \vskip4mm
 \par
   By virtue of (F0) and (F1), for any given $\varepsilon>0$, there exists $C_{\varepsilon}>0$
 such that
 \begin{equation}\label{f20}
   |f(x, t)|\le \varepsilon|t|+C_{\varepsilon}|t|^{p-1}, \ \ \ \ \forall \ (x, t)\in \R^N\times \R.
 \end{equation}
 Consequently,
 \begin{equation}\label{F30}
   |F(x, t)|\le \varepsilon|t|^2+C_{\varepsilon}|t|^{p}, \ \ \ \ \forall \ (x, t)\in \R^N\times \R.
 \end{equation}
 According to \x{F30}, we can demonstrate that $\Phi$ is of class $C^{1}(E, \R)$, and
 \begin{equation}\label{Phd}
   \langle \Phi'(u), v \rangle  =  \int_{{\R}^N}\left(\nabla u\nabla v+V(x)uv\right)\mathrm{d}x
     -\int_{{\R}^N}f(x, u) v\mathrm{d}x, \ \ \ \ \forall \ u, v\in E.
 \end{equation}

 \par
   In the recent paper \cite{Sz}, Szulkin and Weth developed an ingenious approach to find ground state solutions
 for problem \x{SP}. This approach transforms, by  a direct and simple reduction,  the indefinite variational
 problem to a definite one, resulting in a new minimax characterization of the corresponding critical value.
 More precisely, they proved the following theorem.

 \vskip4mm
 \par\noindent
 {\bf Theorem 1.1.} (\cite {Sz})\ \ {\it Assume that $V$ and $f$ satisfy} (V0), (F0), (F1) {\it and the following assumptions:}

 \vskip2mm
 \noindent
 \begin{itemize}
  \item[(V1)]  {\it $V(x)$ is 1-periodic in each of $x_1, x_2, \ldots, x_N$;}

 \item[(F2)]  {\it $f(x, t)$ is 1-periodic in each of $x_1, x_2, \ldots, x_N$;}

 \item[(SQ)] {\it $\lim_{|t|\to \infty}\frac{|F(x, t)|}{|t|^2}=\infty$ uniformly in $x\in \R^{N}$;}

 \item[(Ne)] {\it $t\mapsto \frac{f(x, t)}{|t|}$ is strictly increasing on $(-\infty, 0)\cup (0, \infty)$.}
 \end{itemize}

 \vskip2mm
 \noindent
 {\it Then problem \x{SP} has a solution $u_0\in E$ such that $\Phi(u_0)=m:=\inf_{\mathcal{N}^{-}}\Phi>0$, where}

 \begin{equation}\label{Ne-}
   \mathcal{N}^{-}  = \left\{u\in E\setminus E^{-} : \langle \Phi'(u), u \rangle=\langle \Phi'(u), v \rangle=0,
       \ \forall \ v\in E^{-} \right\}.
 \end{equation}

 \vskip4mm
 \par
   The set $\mathcal{N}^{-}$ was first introduced by Pankov \cite{Pa}, which is a subset of the Nehari manifold
 \begin{equation}\label{Ne}
   \mathcal{N}  = \left\{u\in E\setminus \{0\} : \langle \Phi'(u), u \rangle=0 \right\}.
 \end{equation}
 Since $u_0$ is a solution to the equation $\Phi'(u)=0$ at which $\Phi$ has minimal ``energy" in set $\mathcal{N}^{-}$, we
 shall call it a ground state solution of Nehari-Pankov type. Theorem 1.1 is also established in Pankov
 \cite[Section 5]{Pa} under the following additional assumptions on the nonlinearity: $f\in C^{1}(\R^{N+1}, \R)$,
 $|f'_t(x, t)|\le a_0(1+|t|^{p-2})$ and
 \begin{equation}\label{fd}
   0<\frac{f(x, t)}{t}<\theta f'_t(x, t) \ \ \ \ \mbox{for some}\ \theta\in (0, 1)\ \mbox{and all}\ t\ne 0.
 \end{equation}
 It is easy to see that \x{fd} is stronger than both (Ne) and the following classical condition (AR) due to
 Ambrosetti and Rabinowitz \cite{AR}:

 \vskip2mm
 \begin{itemize}
 \item[(AR)] there exists $\mu>2$ such that
 $$
   0<\mu F(x, t)\le tf(x, t), \ \ \ \ \forall \ x\in \R^N, \ \ t\ne 0.
 $$
 \end{itemize}
 The existence of a nontrivial solution of \x{SP} has been obtained in \cite{BW, CZR, KS, Ra, TW} under (AR)
 and some other standard assumptions of $f$. It is well known that (AR) implies (SQ). The idea of using the more
 natural super-quadratic condition (SQ) to replace (AR) under a Nehari type setting goes back to Liu and Wang \cite{LW}.
 Afterwards, condition (SQ) was also used in many papers, see \cite{CM, DS, DL, HL, JZ, LT, LT1, LWZ, Ta, Ta3, ZZ, ZZW}. In the definite case where
 $\sigma(-\triangle +V) \subset (0, \infty)$,  Theorem 1.1 is a slight extension of a result by Li, Wang and Zeng \cite{LWZ}.

 \par
   There  have been a few new works on the existence of ``ground state solutions" for problem \x{SP} after Szulkin and
 Weth \cite{Sz} in which, various conditions better than (Ne) are obtained, see \cite{Liu, Sc1, Ta1, Ya}.
 However, the ``ground state solutions" for problem \x{SP} in \cite{Liu, Sc1, Ta1, Ya} are in fact a nontrivial solution
 $u_0$ which satisfies $\Phi(u_0)=\inf_{\mathcal{M}}\Phi$, where
 \begin{equation}\label{M}
   \mathcal{M}  = \left\{u\in E\setminus \{0\} : \Phi'(u)=0 \right\}
 \end{equation}
 is a very small subset of $\mathcal{N}^{-}$. In general, it is much more difficult to find a solution
 $u_0$ for \x{SP} which satisfies $\Phi(u_0)=\inf_{\mathcal{N}^{-}}\Phi$ than one satisfying
 $\Phi(u_0)=\inf_{\mathcal{M}}\Phi$.

 \par
    We point out that the Nehari type assumption (Ne) is very crucial in Szulkin and Weth \cite{Sz}.
 In fact, the starting point of their approach is to show that for each $u\in E\setminus E^{-}$, the
 Nehari-Pankov manifold $\mathcal{N}^{-}$ intersects $\hat{E}(u)$ in exactly one point $\hat{m}(u)$, where
 \begin{equation}\label{Eu}
 \hat{E}(u): = E^{-}\oplus \R^{+} u=E^{-}\oplus \R^{+} u^{+}, \ \ \ \ \mbox{and}\ \ \R^{+} = [0,\infty).
 \end{equation}
 The uniqueness of $\hat{m}(u)$ enables one to define a map $u\mapsto \hat{m}(u)$, which is important in
 the remaining proof. If $t\mapsto f(x, t)/|t|$  {\it is not strictly increasing}, then $\hat{m}(u)$ may
 not be unique and their arguments become invalid. This paper intends to address this problem caused by
 the dropping of this ``strictly increasing" condition on $f$. Motivated by the works \cite{LWZ, Pa, Ra, Sz, Ta2, Ta4},
 we will develop a more direct approach to generalize and improve Theorem 1.1 by relaxing assumption
 (Ne) in two cases: the periodic case and the asymptotically periodic case. Unlike the Nahari manifold method,
 our approach lies on finding a minimizing Cerami sequence for $\Phi$ outside $\mathcal{N}^{-}$ by using the
 diagonal method, see Lemma 2.10 and \cite[Lemma 3.8]{Ta2}.

 \vskip2mm
 \par
   Before presenting our theorems, in addition to (V0), (V1), (F0), (F1) and (F2), we introduce the
 following assumptions.

 \vskip2mm
 \noindent
 \begin{itemize}
 \item[(F3)]$\lim_{|t|\to \infty}\frac{|F(x, t)|}{|t|^2}=\infty, \ \ a.e. \ x\in \R^N$;

 \item[(WN)]$t\mapsto \frac{f(x, t)}{|t|}$ is non-decreasing on $(-\infty, 0)\cup (0, \infty)$.
 \end{itemize}

 \vskip4mm
 \par
   We are now in a position  to state the first result of this paper.

 \vskip4mm
 \par\noindent
 {\bf Theorem 1.2.}\ \ {\it Assume that $V$ and $f$ satisfy} (V0), (V1), (F0), (F1), (F2), (F3) {\it and}
 (WN). {\it Then problem \x{SP} has a solution $u_0\in E$ such that $\Phi(u_0)=\inf_{\mathcal{N}^{-}}\Phi>0$.}

 \vskip4mm
 \par
    Next, we assume that (V0$'$) is satisfied, i.e. $V(x)$ is asymptotically periodic. In this case,
 the functional $\Phi$ loses the $\Z^N$-translation invariance, and a powerful equation
 $u^{+}(\cdot+k)=[u(\cdot+k)]^{+}$ for any $k\in \Z^N$ is no longer valid. For the above reasons,
 many effective methods for periodic problems cannot be applied to asymptotically periodic ones. To
 the best of our knowledge, there are no results on the existence of ground state solutions for \x{SP}
 when $V(x)$ is asymptotically periodic. In this paper, we find new tricks to overcome the difficulties
 caused by the dropping of periodicity of $V(x)$.

 \vskip2mm
 \par
   Let $\mathcal{A}_0=-\triangle +V_0$ with spectral family $\{\mathcal{F}(\lambda): -\infty \le \lambda \le +\infty\}$,
 $|\mathcal{A}_0|$ be the absolute value of $\mathcal{A}_0$, and $|\mathcal{A}_0|^{1/2}$ be the
 square root of $|\mathcal{A}_0|$. Let
 \begin{equation}\label{E0}
   E=\mathfrak{D}(|\mathcal{A}_0|^{1/2}), \ \ \ \ E^{\mathcal{F}-}=\mathcal{F}(0-)E, \ \ \ \ E^{\mathcal{F}+}=[id-\mathcal{F}(0)]E.
 \end{equation}
 For any $u\in E$, it is easy to see that $u=u^{\mathcal{F}-}+u^{\mathcal{F}+}$ and
 \begin{equation}\label{A0u}
   \mathcal{A}_0u^{\mathcal{F}-}=-|\mathcal{A}_0|u^{\mathcal{F}-}, \ \ \ \ \mathcal{A}_0u^{\mathcal{F}+}=|\mathcal{A}_0|u^{\mathcal{F}+},
       \ \ \ \ \forall \ u\in E\cap \mathfrak{D}(\mathcal{A}_0),
 \end{equation}
 where
 \begin{equation}\label{u0-}
   u^{\mathcal{F}-}:=\mathcal{F}(0-)u\in E^{\mathcal{F}-}, \ \ \ \  u^{\mathcal{F}+}:=[id-\mathcal{F}(0)]u\in E^{\mathcal{F}+}.
 \end{equation}
 Obviously, If $V_1=0$, i.e. $V=V_0$, then $u^{-}=u^{\mathcal{F}-}$, $u^{+}=u^{\mathcal{F}+}$,
 $E^{-}=E^{\mathcal{F}-}$ and $E^{+}=E^{\mathcal{F}+}$.

 \vskip4mm
 \par
    Under assumption (V0$'$), we can define a new inner product
 \begin{equation} \label{IP0}
    (u, v)_0=\left(|\mathcal{A}_0|^{1/2}u, |\mathcal{A}_0|^{1/2}v\right)_{L^2}, \ \ \ \ \forall \ u, v\in E
 \end{equation}
 and the corresponding norm
 \begin{equation} \label{N0}
    \|u\|_0=\left\||\mathcal{A}_0|^{1/2}u\right\|_{2}, \ \ \ \ \forall \ u\in E.
 \end{equation}
 It is easy to see the norm $\|\cdot\|_0$ is equivalent to the norm $\|\cdot\|_{H^1(\R^N)}$. In particular,
 by (V0$'$), one has
 \begin{equation} \label{N+0}
    \|u\|_0^2\ge \Pi_0\|u\|_2^2, \ \ \ \ \forall \ u\in E^{\mathcal{F}+},
 \end{equation}
 where
 \begin{equation} \label{Pi}
  \Pi_0 := \inf[\sigma(-\triangle +V_0)\cap (0, \infty)].
 \end{equation}

 \vskip4mm
 \par
   Instead of (V1) and (F2), we make the following assumptions.

 \vskip2mm
 \noindent
 \begin{itemize}
 \item[(V1$'$)] $V_0(x)$ is 1-periodic in each of $x_1, x_2, \ldots, x_N$, and
 $$
   0\le -V_1(x)\le \sup_{\R^N}[-V_1(x)]<\Pi_0, \ \ \ \ \forall \ x\in \R^N;
 $$

 \item[(F2$'$)] $f(x, t)=f_0(x, t)+f_1(x, t)$, $f_0\in C(\R^N\times \R)$, $f_0(x, t)$ is 1-periodic in each of
 $x_1, x_2, \ldots, x_N$, $t\mapsto f_0(x, t)/|t|$ is non-decreasing on $(-\infty, 0)\cup (0, \infty)$; $f_1\in C(\R^N\times \R)$
 satisfies that
 $$
   0\le tf_1(x, t)\le a(x)\left(|t|^2+|t|^{p}\right), \ \ \ \ \forall \ (x, t)\in \R^N\times \R,
 $$
 where $a\in C({\R}^{N})$ with $\lim_{|x|\to\infty}a(x)=0$.
 \end{itemize}

 \vskip4mm
 \par
   We are now in a position to state the second result of this paper.

 \vskip4mm
 \par\noindent
 {\bf Theorem 1.3.}\ \ {\it Assume that $V$ and $f$ satisfy} (V0$'$), (V1$'$), (F0), (F1), (F2$'$), (SQ) {\it and}
 (WN). {\it Moreover, $-V_1(x)t^2+F_1(x, t)>0$ for $|x|<1+\sqrt{N}$ and $t\ne 0$. Then problem \x{SP} has a
 solution $u_0\in E$ such that $\Phi(u_0)=\inf_{\mathcal{N}^{0}}\Phi>0$, where}
 \begin{equation}\label{Ne0}
   \mathcal{N}^{0}  = \left\{u\in E\setminus E^{\mathcal{F}-} : \langle \Phi'(u), u \rangle=\langle \Phi'(u), v \rangle=0,
       \ \forall \ v\in E^{\mathcal{F}-} \right\}.
 \end{equation}

 \vskip4mm
 \par
   The remainder of this paper is organized as follows. In Section 2, some preliminary results are presented.
 The proofs of Theorems 1.2 and 1.3 are given in Section 3 and Section 4, respectively.


 \vskip6mm
 {\section{Preliminaries}}
 \setcounter{equation}{0}

 \par
   Let $X$ be a real Hilbert space with $X=X^{-}\oplus X^{+}$ and $X^{-}\bot\  X^{+}$.
 For a functional $\varphi\in C^{1}(X, \R)$, $\varphi$ is said to be weakly sequentially
 lower semi-continuous if for any $u_n\rightharpoonup u$ in $X$ one has $\varphi(u)\le \liminf_{n\to\infty}\varphi(u_n)$,
 and $\varphi'$ is said to be weakly sequentially continuous if $\lim_{n\to\infty}\langle\varphi'(u_n), v\rangle=
 \langle\varphi'(u), v\rangle$ for each $v\in X$.

 \vskip4mm
 \par\noindent
 {\bf Lemma 2.1.} (\cite {KS, LS})\ \ {\it Let $X$ be a real Hilbert space with $X=X^{-}\oplus X^{+}$ and $X^{-}\bot\  X^{+}$,
 and let $\varphi\in C^{1}(X, \R)$ of the form
 $$
   \varphi(u)=\frac{1}{2}\left(\|u^{+}\|^2-\|u^{-}\|^2\right)-\psi(u), \ \ \ \ u=u^{-}+u^{+}\in X^{-}\oplus X^{+}.
 $$
 Suppose that the following assumptions are satisfied:}
 \vskip2mm
 \par
 (KS1)\ \ {\it $\psi\in C^{1}(X, \R)$ is bounded from below and weakly sequentially lower semi-continuous;}

 \vskip2mm
 \par
 (KS2)\ \ {\it $\psi'$ is weakly sequentially continuous;}

 \vskip2mm
 \par
 (KS3)\ \ {\it there exist $r>\rho>0$ and $e\in X^{+}$ with $\|e\|=1$ such that
 $$
   \kappa:=\inf\varphi(S^{+}_{\rho}) > \sup \varphi(\partial Q),
 $$
 where }
 $$
     S^{+}_{\rho}=\left\{u\in X^{+} : \|u\|=\rho\right\},  \ \ \ \ Q=\left\{v+se : v\in X^{-},\ s\ge 0,\ \|v+se\|\le r\right\}.
 $$

 \vskip2mm
 \noindent
 {\it Then there exist a constant $c\in [\kappa, \sup \varphi(Q)]$ and a sequence $\{u_n\}\subset X$ satisfying}
 \begin{eqnarray*}
   \varphi(u_n)\rightarrow c, \ \ \ \ \|\varphi'(u_n)\|(1+\|u_n\|)\rightarrow 0.
 \end{eqnarray*}

 \vskip4mm
 \par
    We set
 \begin{equation}\label{Psi}
   \Psi(u)=\int_{{\R}^N}F(x, u)\mathrm{d}x, \ \ \ \ \forall \ u\in E.
 \end{equation}

 \vskip4mm
 \par
   Employing a standard argument, one checks easily the following:

 \vskip4mm
 \par\noindent
 {\bf Lemma 2.2.}\ \ {\it Suppose that} (V0$'$), (F0) {\it and} (F1) {\it  are satisfied, and $F(x, t)\ge 0$ for
 all $(x, t)\in \R^N\times \R$. Then $\Psi$ is nonnegative, weakly sequentially lower semi-continuous, and
 $\Psi'$ is weakly sequentially continuous.}

 \vskip4mm
 \par\noindent
 {\bf Lemma 2.3.}\ \ {\it Suppose that} (V0$'$), (F0), (F1) {\it and} (WN) {\it are satisfied. Then}
 \begin{eqnarray}\label{L31}
   \Phi(u) & \ge & \Phi(tu+w)+\frac{1}{2}\|w\|_0^2-\frac{1}{2}\int_{{\R}^N}V_1(x)w^2\mathrm{d}x\nonumber\\
           &     & \ \ +\frac{1-t^2}{2}\langle\Phi'(u), u \rangle-t\langle\Phi'(u), w \rangle,
                   \ \ \ \ \forall \ u\in E, \ \ t\ge 0, \ \ w\in E^{\mathcal{F}-}.
 \end{eqnarray}

 \vskip2mm
 \par\noindent
 {\bf Proof.} \ \ For any $x\in \R^N$ and $\tau\ne 0$, (WN) yields
 \begin{equation}\label{L30}
   f(x, s)\le \frac{f(x, \tau)}{|\tau|}|s|, \ \ \ \  s \le \tau; \ \ \ \ f(x, s)\ge \frac{f(x, \tau)}{|\tau|}|s|, \ \ \ \ s\ge \tau.
 \end{equation}
 It follows that
 \begin{equation}\label{L32}
   \left(\frac{1-t^2}{2}\tau-t\sigma \right)f(x, \tau)
      \ge \int_{t\tau+\sigma}^{\tau}f(x, s)\mathrm{d}s,  \ \ \ \ t\ge 0, \ \ \sigma\in \R.
 \end{equation}
 To show \x{L32}, we consider four possible cases. By virtue of \x{L30} and $sf(x, s)\ge 0$, one has

 \vskip2mm
 \par
   Case 1). \ $0\le t\tau+\sigma\le \tau$ or $t\tau+\sigma\le \tau \le 0$,
 $$
   \int_{t\tau+\sigma}^{\tau}f(x, s)\mathrm{d}s\le \frac{f(x, \tau)}{|\tau|}\int_{t\tau+\sigma}^{\tau}|s|\mathrm{d}s
     \le \left(\frac{1-t^2}{2}\tau-t\sigma \right)f(x, \tau);
 $$

 \vskip2mm
 \par
   Case 2). \ $t\tau+\sigma\le 0 \le \tau$,
 $$
   \int_{t\tau+\sigma}^{\tau}f(x, s)\mathrm{d}s\le \int_{0}^{\tau}f(x, s)\mathrm{d}s
     \le \frac{f(x, \tau)}{|\tau|}\int_{0}^{\tau}|s|\mathrm{d}s
     \le \left(\frac{1-t^2}{2}\tau-t\sigma \right)f(x, \tau);
 $$

 \vskip2mm
 \par
   Case 3). \ $0 \le \tau\le t\tau+\sigma$ or $\tau\le t\tau+\sigma\le 0$,
 $$
   \int_{\tau}^{t\tau+\sigma}f(x, s)\mathrm{d}s\ge \frac{f(x, \tau)}{|\tau|}\int_{\tau}^{t\tau+\sigma}|s|\mathrm{d}s
     \ge -\left(\frac{1-t^2}{2}\tau-t\sigma \right)f(x, \tau);
 $$

 \vskip2mm
 \par
   Case 4). \ $\tau \le 0 \le t\tau+\sigma$,
 $$
   \int_{\tau}^{t\tau+\sigma}f(x, s)\mathrm{d}s\ge \int_{\tau}^{0}f(x, s)\mathrm{d}s
     \ge \frac{f(x, \tau)}{|\tau|}\int_{\tau}^{0}|s|\mathrm{d}s
     \ge -\left(\frac{1-t^2}{2}\tau-t\sigma \right)f(x, \tau).
 $$
 The above four cases show that \x{L32} holds.

 \par
   We let $b: E\times E \rightarrow \R$ denote the symmetric bilinear form given by
 \begin{equation}\label{buv}
   b(u, v)=\int_{{\R}^N}(\nabla u\nabla v+V(x)uv)\mathrm{d}x, \ \ \ \ \forall \ u, v\in E.
 \end{equation}
 By virtue of \x{Ph}, \x{Phd} and \x{buv}, one has
 \begin{equation}\label{Ph2}
   \Phi(u) = \frac{1}{2}b(u, u)-\int_{{\R}^N}F(x, u)\mathrm{d}x, \ \ \ \ \forall \ u\in E
 \end{equation}
 and
 \begin{equation}\label{Phd2}
   \langle \Phi'(u), v \rangle  = b(u, v)-\int_{{\R}^N}f(x, u)v\mathrm{d}x, \ \ \ \ \forall \ u, v\in E.
 \end{equation}
 Thus, by \x{A0u}, \x{N0}, \x{L32}, \x{Ph2} and \x{Phd2}, one has
 \begin{eqnarray*}
   &     & \Phi(u)-\Phi(tu+w)\\
   &  =  & \frac{1}{2}[ b(u, u)-b(tu+w, tu+w)]+\int_{\R^N}[F(x, tu+w)-F(x, u)]\mathrm{d}x\\
   &  =  & \frac{1-t^2}{2} b(u, u)-tb(u, w)-\frac{1}{2}b(w, w)+\int_{\R^N}[F(x, tu+w)-F(x, u)]\mathrm{d}x\\
   &  =  & -\frac{1}{2}b(w, w)+\frac{1-t^2}{2}\langle\Phi'(u), u \rangle-t\langle\Phi'(u), w \rangle\\
   &     &   +\int_{\R^N}\left[\frac{1-t^2}{2}f(x, u)u-tf(x, u)w-\int_{tu+w}^{u}f(x, s)\mathrm{d}s\right]\mathrm{d}x\\
   &  =  & \frac{1}{2}\|w\|_0^2-\frac{1}{2}\int_{{\R}^N}V_1(x)w^2\mathrm{d}x
             +\frac{1-t^2}{2}\langle\Phi'(u), u \rangle-t\langle\Phi'(u), w \rangle\\
   &     &   +\int_{\R^N}\left[\frac{1-t^2}{2}f(x, u)u-tf(x, u)w-\int_{tu+w}^{u}f(x, s)\mathrm{d}s\right]\mathrm{d}x\\
   & \ge & \frac{1}{2}\|w\|_0^2-\frac{1}{2}\int_{{\R}^N}V_1(x)w^2\mathrm{d}x+\frac{1-t^2}{2}\langle\Phi'(u), u \rangle-t\langle\Phi'(u), w \rangle,
            \ \ \ \ \forall \ t\ge 0, \ \ w\in E^{\mathcal{F}-}.
 \end{eqnarray*}
 This shows that \x{L31} holds.
 \qed

 \vskip4mm
 \par
    From Lemma 2.3, we have the following two corollaries.

 \vskip4mm
 \par\noindent
 {\bf Corollary 2.4.}\ \ {\it Suppose that} (V0$'$), (F0), (F1) {\it and} (WN) {\it are satisfied. Then for
 $u\in \mathcal{N}^{0}$}
 \begin{equation}\label{2701}
   \Phi(u) \ge \Phi(tu+w)+\frac{1}{2}\|w\|_0^2-\frac{1}{2}\int_{{\R}^N}V_1(x)w^2\mathrm{d}x,
       \ \ \ \ \forall \ t\ge 0, \ \ w\in E^{\mathcal{F}-}.
 \end{equation}

 \vskip4mm
 \par\noindent
 {\bf Corollary 2.5.}\ \ {\it Suppose that} (V0$'$), (F0), (F1) {\it and} (WN) {\it are satisfied. Then}
 \begin{eqnarray}\label{2801}
   \Phi(u) & \ge & \frac{t^2}{2}\|u\|_0^2-\int_{{\R}^N}F(x, tu^{\mathcal{F}+})\mathrm{d}x
                     +\frac{1-t^2}{2}\langle\Phi'(u), u \rangle+t^2\langle\Phi'(u), u^{\mathcal{F}-} \rangle\nonumber\\
           &     & \ \  +\frac{t^2}{2}\int_{{\R}^N}V_1(x)\left[(u^{\mathcal{F}+})^2-(u^{\mathcal{F}-})^2\right]\mathrm{d}x,
                     \ \ \ \ \forall \ u\in E, \ \ t\ge 0. \ \ \ \
 \end{eqnarray}

 \vskip4mm
 \par\noindent
 {\bf Lemma 2.6.} \ \ {\it Suppose that} (V0$'$), (V1$'$), (F0), (F1) {\it and} (WN) {\it are satisfied. Then}

 \vskip2mm
 \par
  \ \ (i) \ {\it there exists $\rho>0$ such that }
 \begin{equation}\label{m1}
   m:=\inf_{\mathcal{N}^{0}}\Phi \ge \kappa:=\inf \left\{\Phi(u) : u\in E^{\mathcal{F}+}, \|u\|_0=\rho\right\}>0.
 \end{equation}

 \vskip2mm
 \par
  \ \ (ii) \ {\it $\|u^{\mathcal{F}+}\|_0\ge \max\left\{\|u^{\mathcal{F}-}\|_0, \sqrt{2m}\right\}$ for all $u\in \mathcal{N}^{0}$.}

 \vskip2mm
 \par\noindent
 {\bf Proof.} \ \ (i) \ Set $\Theta_0=\sup_{\R^N}[-V_1(x)]$. Let $\varepsilon_0=(\Pi_0-\Theta_0)/3$.
 Then \x{F30} implies that there exists a constant $C_{\varepsilon_0}>0$ such that
 \begin{equation}\label{2901}
   F(x, tu^{\mathcal{F}+})\le \varepsilon_0|tu^{\mathcal{F}+}|^2+C_{\varepsilon_0}|tu^{\mathcal{F}+}|^{p}, \ \ \ \ \forall \ x\in \R^N.
 \end{equation}
 From \x{Ph}, \x{N+0}, \x{2701} and \x{2901}, we have for $u\in \mathcal{N}^{0}$
 \begin{eqnarray}\label{2902}
   \Phi(u)
     & \ge & \Phi(tu^{\mathcal{F}+})\nonumber\\
     &  =  & \frac{t^2}{2}\|u^{\mathcal{F}+}\|_0^2+\frac{t^2}{2}\int_{{\R}^N}V_1(x)(u^{\mathcal{F}+})^2\mathrm{d}x
               -\int_{\R^N}F(x, tu^{\mathcal{F}+})\mathrm{d}x\nonumber\\
     & \ge & \frac{t^2}{2}\left[\|u^{\mathcal{F}+}\|_0^2-(\Theta_0+2\varepsilon_0)\|u^{\mathcal{F}+}\|_2^2\right]
               -t^pC_{\varepsilon_0}\|u^{\mathcal{F}+}\|_p^p\nonumber\\
     & \ge & \frac{t^2}{2}\left(1-\frac{\Theta_0+2\varepsilon_0}{\Pi_0}\right)\|u^{\mathcal{F}+}\|_0^2
               -t^p\gamma_p^pC_{\varepsilon_0}\|u^{\mathcal{F}+}\|_0^p>0 \ \ \ \ \mbox{for small} \ t>0.
 \end{eqnarray}
 This shows that there exists a $\rho>0$ such that \x{m1} holds.

 \par
    (ii) \ (F1) and (WN) imply that $F(x, u)\ge 0$. Hence, it follows from (i), \x{Ph} and the definition of $\mathcal{N}^{0}$ that (ii) holds.
 \qed

 \vskip4mm
 \par\noindent
 {\bf Lemma 2.7.}\ \ {\it Suppose that} (V0$'$), (V1$'$), (F0), (F1) {\it and} (F3) {\it are satisfied. Then for any
 $e\in E^{\mathcal{F}+}$, $\sup \Phi(E^{\mathcal{F}-}\oplus \R^{+} e)<\infty$, and there is $R_e>0$ such that}
 \begin{equation}\label{L321}
   \Phi(u) \le 0, \ \ \ \ \forall \   u\in E^{\mathcal{F}-}\oplus \R^{+} e, \ \ \|u\|_0\ge R_e.
 \end{equation}

 \vskip2mm
 \par\noindent
 {\bf Proof.} \ \ Arguing indirectly, assume that for some sequence $\{w_n+s_ne\}\subset
 E^{\mathcal{F}-}\oplus \R^{+} e$ with $\|w_n+s_ne\|_0\rightarrow \infty$ such that $\Phi(w_n+s_ne) \ge 0$ for all $n\in \N$.
 Set $v_n=(w_n+s_ne)/\|w_n+s_ne\|_0=v_n^{\mathcal{F}-}+t_ne$, then $\|v_n^{\mathcal{F}-}+t_ne\|_0=1$.
 Passing to a subsequence, we may assume that $t_n\rightarrow \bar{t}$, $v_n^{\mathcal{F}-}\rightharpoonup v^{\mathcal{F}-}$,
 and $v_n^{\mathcal{F}-}\rightarrow v^{\mathcal{F}-}$ a.e. on $\R^N$. Hence,
 \begin{eqnarray}\label{L402}
   0 & \le & \frac{\Phi(w_n+s_ne)}{\|w_n+s_nu\|_0^2}\nonumber\\
     &  =  & \frac{t_n^2}{2}\|e\|_0^2-\frac{1}{2}\|v_n^{\mathcal{F}-}\|_0^2+\frac{1}{2}\int_{{\R}^N}V_1(x)(v_n^{\mathcal{F}-}+t_ne)^2\mathrm{d}x
              -\int_{\R^N}\frac{F(x, w_n+s_ne)}{\|w_n+s_ne\|_0^2}\mathrm{d}x\nonumber\\
     & \le & \frac{t_n^2}{2}\|e\|_0^2-\frac{1}{2}\|v_n^{\mathcal{F}-}\|_0^2
              -\int_{\R^N}\frac{F(x, w_n+s_ne)}{\|w_n+s_ne\|_0^2}\mathrm{d}x.
 \end{eqnarray}

 \par
   If $\bar{t}=0$, then it follows from \x{L402} that
 $$
   0 \le \frac{1}{2}\|v_n^{\mathcal{F}-}\|_0^2+\int_{\R^N}\frac{F(x, w_n+s_ne)}{\|w_n+s_ne\|_0^2}\mathrm{d}x
     \le \frac{t_n^2}{2}\|e\|_0^2 \rightarrow 0,
 $$
 which yields $\|v_n^{\mathcal{F}-}\|_0\rightarrow 0$, and so $1=\|v_n^{\mathcal{F}-}+t_ne\|_0^2 \rightarrow 0$, a contradiction.

 \par
   If $\bar{t}\ne 0$, then it follows from \x{L402} and (F3) that
 \begin{eqnarray*}
  0 & \le & \limsup_{n\to\infty}\left[\frac{t_n^2}{2}\|e\|_0^2-\frac{1}{2}\|v_n^{\mathcal{F}-}\|_0^2
              -\int_{\R^N}\frac{F(x, w_n+s_ne)}{\|w_n+s_ne\|_0^2}\mathrm{d}x\right]\\
    & \le & \frac{\bar{t}^2}{2}\|e\|_0^2
             -\liminf_{n\to\infty}\int_{\R^N}\frac{F(x, w_n+s_ne)}{(w_n+s_ne)^2}\left(v_n^{\mathcal{F}-}+t_ne\right)^2\mathrm{d}x\\
    & \le & \frac{\bar{t}^2}{2}\|e\|_0^2
             -\int_{\R^N}\liminf_{n\to\infty}\frac{F(x, w_n+s_ne)}{(w_n+s_ne)^2}\left(v_n^{\mathcal{F}-}+t_ne\right)^2\mathrm{d}x\\
    &  =  & -\infty,
 \end{eqnarray*}
 a contradiction.
 \qed

 \vskip4mm
 \par\noindent
 {\bf Corollary 2.8.}\ \ {\it Suppose that} (V0$'$), (V1$'$), (F0), (F1) {\it and} (F3) {\it are satisfied.
 Let $e\in E^{\mathcal{F}+}$ with $\|e\|_0=1$. Then there is a $r_0>\rho$ such that $\sup \Phi(\partial Q)\le 0$
 for $r\ge r_0$, where}
 \begin{equation}\label{Q1}
   Q=\left\{w+se : w\in E^{\mathcal{F}-}, \ s\ge 0, \|w+se\|_0\le r\right\}.
 \end{equation}

 \vskip4mm
 \par\noindent
 {\bf Lemma 2.9.}\ \ {\it Suppose that} (V0$'$), (V1$'$), (F0), (F1), (F3) {\it and} (WN) {\it are satisfied.
 Then there exist a constant $c\in [\kappa, \sup \Phi(Q)]$ and a sequence $\{u_n\}\subset E$ satisfying
 \begin{equation}\label{Ce1}
   \Phi(u_n)\rightarrow c, \ \ \ \ \|\Phi'(u_n)\|(1+\|u_n\|_0)\rightarrow 0,
 \end{equation}
 where $Q$ is defined by \x{Q1}.}

 \vskip2mm
 \par\noindent
 {\bf Proof.} \ \ Lemma 2.9 is a direct corollary of Lemmas 2.1, 2.2, 2.6 (i) and Corollary 2.8.
 \qed

 \vskip4mm
 \par
    The following lemma is crucial to demonstrate the existence of ground state solutions of Nehari-Pankov type for problem \x{SP}.

 \vskip4mm
 \par\noindent
 {\bf Lemma 2.10.}\ \ {\it Suppose that} (V0$'$), (V1$'$), (F0), (F1), (F3) {\it and} (WN) {\it are satisfied.
 Then there exist a constant $c_*\in [\kappa, m]$ and a sequence $\{u_n\}\subset E$ satisfying}
 \begin{equation}\label{Ce}
   \Phi(u_n)\rightarrow c_*, \ \ \ \ \|\Phi'(u_n)\|(1+\|u_n\|_0)\rightarrow 0.
 \end{equation}

 \vskip2mm
 \par\noindent
 {\bf Proof.} \ \ Choose $v_k\in \mathcal{N}^{0}$ such that
 \begin{equation}\label{1201}
   m\le \Phi(v_k)< m+\frac{1}{k}, \ \ \ \ k\in \N.
 \end{equation}
 By Lemma 2.6, $\|v_k^{\mathcal{F}+}\|_0\ge \sqrt{2m}>0$. Set $e_k=v_k^{\mathcal{F}+}/\|v_k^{\mathcal{F}+}\|_0$.
 Then $e_k\in E^{\mathcal{F}+}$ and $\|e_k\|_0=1$. In view of Corollary 2.8, there exists $r_k>\max\{\rho, \|v_k\|_0\}$
 such that $\sup \Phi(\partial Q_k)\le 0$, where
 \begin{equation}\label{1202}
   Q_k=\{w+se_k : w\in E^{\mathcal{F}-}, \  s\ge 0,\ \|w+se_k\|_0\le r_k\}, \ \ \ \ k\in \N.
 \end{equation}
 Hence, applying Lemma 2.9 to the above set $Q_k$, there exist a constant $c_k\in [\kappa, \sup \Phi(Q_k)]$ and
 a sequence $\{u_{k, n}\}_{n\in \N}\subset E$ satisfying
 \begin{equation}\label{Cek}
   \Phi(u_{k, n})\rightarrow c_k, \ \ \ \ \|\Phi'(u_{k, n})\|(1+\|u_{k, n}\|_0)\rightarrow 0, \ \ \ \ k\in \N.
 \end{equation}
 By virtue of Corollary 2.4, one can get that
 \begin{equation}\label{1203}
   \Phi(v_k) \ge \Phi(tv_k+w), \ \ \ \ \forall \ t\ge 0, \ \ w\in E^{\mathcal{F}-}.
 \end{equation}
 Since $v_k\in Q_k$, it follows from \x{1202} and \x{1203} that $\Phi(v_k)=\sup \Phi(Q_k)$. Hence, by \x{1201}
 and \x{Cek}, one has
 \begin{equation}\label{1204}
   \Phi(u_{k, n})\rightarrow c_k< m+\frac{1}{k}, \ \ \ \ \|\Phi'(u_{k, n})\|(1+\|u_{k, n}\|_0)\rightarrow 0, \ \ \ \ k\in \N.
 \end{equation}
 Now, we can choose a sequence $\{n_k\}\subset \N$ such that
 \begin{equation}\label{1205}
   \Phi(u_{k, n_k})< m+\frac{1}{k}, \ \ \ \ \|\Phi'(u_{k, n_k})\|(1+\|u_{k, n_k}\|_0)<\frac{1}{k}, \ \ \ \ k\in \N.
 \end{equation}
 Let $u_k=u_{k, n_k}, k\in \N$. Then, going if necessary to a subsequence, we have
 \begin{equation*}\label{1206}
   \Phi(u_n)\rightarrow c_*\in [\kappa, m], \ \ \ \ \|\Phi'(u_n)\|(1+\|u_n\|_0)\rightarrow 0.
 \end{equation*}
 \qed

 \vskip4mm
 \par\noindent
 {\bf Lemma 2.11.}\ \ {\it Suppose that} (V0$'$), (V1$'$), (F0), (F1), (SQ) {\it and} (WN) {\it are satisfied. Then for any
 $u\in E\setminus E^{\mathcal{F}-}$, there exist $t(u)>0$ and $w(u)\in E^{\mathcal{F}-}$ such that $t(u)u+w(u)\in \mathcal{N}^{0}$.}

 \vskip2mm
 \par\noindent
 {\bf Proof.} \ \ Since $E^{\mathcal{F}-}\oplus\R^{+}u=E^{\mathcal{F}-}\oplus \R^{+}u^{\mathcal{F}+}$, we may assume that
 $u\in E^{\mathcal{F}+}$. By Lemma 2.7, there exists $R>0$ such that $\Phi(u)\le 0$ for $u\in (E^{\mathcal{F}-}
 \oplus \R^{+} u)\setminus B_{R}(0)$. By Lemma 2.6 (i), $\Phi(tu)>0$ for small $t>0$. Thus, $0<\sup\Phi(E^{\mathcal{F}-}
 \oplus \R^{+} u)<\infty$. It is easy see that $\Phi$ is weakly upper semi-continuous
 on $E^{\mathcal{F}-}\oplus \R^{+} u$, therefore, $\Phi(u_0)=\sup\Phi(E^{\mathcal{F}-}\oplus \R^{+} u)$ for some
 $u_0\in E^{\mathcal{F}-}\oplus \R^{+} u$. This $u_0$ is a critical point of $\Phi|_{E^{\mathcal{F}-}\oplus \R u}$, so
 $\langle \Phi'(u_0), u_0 \rangle= \langle \Phi'(u_0), v \rangle=0$ for all $v\in E^{\mathcal{F}-}\oplus \R\ u$.
 Consequently, $u_0\in \mathcal{N}^{0}\cap (E^{\mathcal{F}-}\oplus \R^{+} u)$.
 \qed

 \vskip6mm
 {\section{The periodic case}}
 \setcounter{equation}{0}

 \vskip4mm
 \par
   In this section, we assume that $V$ and $f$ are 1-periodic in each of $x_1, x_2, \ldots, x_N$, i.e., (V1) and
 (F2) are satisfied. In this case, $V_0=V$ and $V_1=0$. Thus, $u^{-}=u^{\mathcal{F}-}$, $u^{+}=u^{\mathcal{F}+}$,
 $E^{-}=E^{\mathcal{F}-}$, $E^{+}=E^{\mathcal{F}+}$ and $\|u\|=\|u\|_0$.

 \vskip4mm
 \par\noindent
 {\bf Lemma 3.1.}\ \ {\it Suppose that} (V0), (V1), (F0), (F1), (F2), (F3) {\it and} (WN) {\it are satisfied. Then any sequence
 $\{u_n\}\subset E$ satisfying
 \begin{equation}\label{Ce2}
   \Phi(u_n)\rightarrow c\ge 0, \ \ \ \ \langle\Phi'(u_n), u_n^{\pm} \rangle\rightarrow 0
 \end{equation}
 is bounded in $E$.}

 \vskip2mm
 \par\noindent
 {\bf Proof.} \ \ To prove the boundedness of $\{u_n\}$, arguing by contradiction, suppose that
 $\|u_n\| \to \infty$. Let $v_n=u_n/\|u_n\|$, then $\|v_n\|=1$. By Sobolev imbedding theorem, there
 exists a constant $C_2>0$ such that $\|v_n\|_2 \le C_2$. If
 $$
   \delta:=\limsup_{n\to\infty}\sup_{y\in \R^N}\int_{B_1(y)}|v_n^{+}|^2\mathrm{d}x=0,
 $$
 then by Lions' concentration compactness principle \cite{Lio} or \cite[Lemma 1.21]{Wi}, $v_n^{+}\rightarrow 0$ in
 $L^{s}(\R^N)$ for $2<s<2^*$. Fix $R>[2(1+c)]^{1/2}$. By virtue of (F0) and (F1), for  $\varepsilon =1/4(RC_2)^2>0$,
 there exists $C_{\varepsilon}>0$ such that \x{F30} holds. Hence, it follows that
 \begin{eqnarray}\label{L61}
   \limsup_{n\to\infty}\int_{\R^N}F(x, Rv_n^{+})\mathrm{d}x
     & \le & \limsup_{n\to\infty}\left[\varepsilon R^2\|v_n^{+}\|_2^2+C_{\varepsilon}R^p\|v_n^{+}\|_p^{p}\right]\nonumber\\
     & \le & \varepsilon(RC_2)^2=\frac{1}{4}.
 \end{eqnarray}
 Let $t_n=R/\|u_n\|$. Hence, by virtue of \x{Ce2}, \x{L61} and Corollary 2.5, one can get that
 \begin{eqnarray*}
   c+o(1)
   &  =  & \Phi(u_n)\\
   & \ge & \frac{t_n^2}{2}\|u_n\|^2-\int_{\R^N}F(x, t_nu_n^{+})\mathrm{d}x
             +\frac{1-t_n^2}{2}\langle\Phi'(u_n), u_n \rangle\\
   &     & \ \  +t_n^2\langle\Phi'(u_n), u_n^{-} \rangle\\
   &  =  & \frac{R^2}{2}\|v_n\|^2-\int_{\R^N}F(x, Rv_n^{+})\mathrm{d}x
             +\left(\frac{1}{2}-\frac{R^2}{2\|u_n\|^2}\right)\langle\Phi'(u_n), u_n \rangle\\
   &     & \ \  +\frac{R^2}{\|u_n\|^2}\langle\Phi'(u_n), u_n^{-} \rangle\\
   &  =  & \frac{R^2}{2}-\int_{\R^N}F(x, Rv_n^{+})\mathrm{d}x +o(1)\\
   & \ge & \frac{R^2}{2}-\frac{1}{4} +o(1) >  c+\frac{3}{4} +o(1).
 \end{eqnarray*}
 This contradiction shows that $\delta>0$.

 \par
    Passing to a subsequence, we may assume the existence of $k_n\in \Z^N$ such that
 $\int_{B_{1+\sqrt{N}}(k_n)}|v_n^{+}|^2\mathrm{d}x> \frac{\delta}{2}$. Let $w_n(x)=v_n(x+k_n)$. Since
 $V(x)$ is 1-periodic in each of $x_1, x_2, \ldots, $ $x_N$. Then
 $\|w_n\|=\|v_n\|=1$, and
 \begin{equation}\label{L62}
   \int_{B_{1+\sqrt{N}}(0)}|w_n^{+}|^2\mathrm{d}x> \frac{\delta}{2}.
 \end{equation}
 Passing to a subsequence, we have $w_n\rightharpoonup w$ in $E$, $w_n\rightarrow w$ in $L^{s}_{\mathrm{loc}}(\R^N)$,
 $2\le s<2^*$, $w_n\rightarrow w$ a.e. on $\R^N$. Thus, \x{L62} implies that $w\ne 0$.

 \par
    Now we define $u_n^{k_n}(x)=u_n(x+k_n)$, then $u_n^{k_n}/\|u_n\|=w_n\rightarrow w$ a.e. on $\R^N$, $w\ne 0$.
 For $x\in \{y\in \R^N : w(y)\ne 0\}$, we have $\lim_{n\to\infty}|u_n^{k_n}(x)|=\infty$. Hence, it follows
 from \x{Ce2}, (F1), (F2), (F3), (WN) and Fatou's lemma that
 \begin{eqnarray*}
  0 &  =  & \lim_{n\to\infty}\frac{c+o(1)}{\|u_n\|^2} = \lim_{n\to\infty}\frac{\Phi(u_n)}{\|u_n\|^2}\\
    &  =  & \lim_{n\to\infty}\left[\frac{1}{2}\left(\|v_n^{+}\|^2-\|v_n^{-}\|^2\right)
             -\int_{\R^N}\frac{F(x, u_n^{k_n})}{(u_n^{k_n})^2}w_n^2\mathrm{d}x\right]\\
    & \le & \frac{1}{2}-\liminf_{n\to\infty}\int_{\R^N}\frac{F(x, u_n^{k_n})}{(u_n^{k_n})^2}w_n^2\mathrm{d}x
             \le \frac{1}{2}-\int_{\R^N}\liminf_{n\to\infty}\frac{F(x, u_n^{k_n})}{(u_n^{k_n})^2}w_n^2\mathrm{d}x  = -\infty.
 \end{eqnarray*}
 This contradiction shows that $\{u_n\}$ is bounded.
 \qed

 \vskip4mm
 \par\noindent
 {\bf Proof of Theorem 1.2.} \  Applying Lemmas 2.10 and 3.1, we deduce that there exists a bounded sequence
 $\{u_n\}\subset E$ satisfying \x{Ce}. Thus there exists a constant $C_3>0$ such that
 $\|u_n\|_2\le C_3$.  If
 $$
   \delta:=\limsup_{n\to\infty}\sup_{y\in \R^N}\int_{B_1(y)}|u_n|^2\mathrm{d}x=0,
 $$
 then by Lions' concentration compactness principle \cite{Lio} or \cite[Lemma 1.21]{Wi}, $u_n\rightarrow 0$ in $L^{s}(\R^N)$ for
 $2<s<2^*$. By virtue of (F0) and (F1), for  $\varepsilon =c_*/4C_3^2>0$ there exists $C_{\varepsilon}>0$ such that \x{f20} and
 \x{F30} hold. It follows that
 \begin{equation}\label{L72}
   \limsup_{n\to\infty}\int_{\R^N}\left[\frac{1}{2}f(x, u_n)u_n-F(x, u_n)\right]\mathrm{d}x\le \frac{3\varepsilon}{2}C_3^2
     +\frac{3}{2}C_{\varepsilon}\lim_{n\to\infty}\|u_n\|_p^{p}=\frac{3c_*}{8}.
 \end{equation}
 From \x{Ph}, \x{Phd}, \x{Ce} and \x{L72}, one can get that
 \begin{eqnarray*}
   c_*
   &  =  & \Phi(u_n)-\frac{1}{2}\langle\Phi'(u_n), u_n \rangle+o(1)\\
   &  =  & \int_{\R^N}\left[\frac{1}{2}f(x, u_n)u_n-F(x, u_n)\right]\mathrm{d}x+o(1)\le \frac{3c_*}{8}+o(1),
 \end{eqnarray*}
 which is a contradiction. Thus $\delta>0$.

 \par
   Passing to a subsequence, we may assume the existence of $k_n\in \Z^N$ such that
 $\int_{B_{1+\sqrt{N}}(k_n)}|u_n|^2\mathrm{d}x> \frac{\delta}{2}$. Let us define $v_n(x)=u_n(x+k_n)$ so that
 \begin{equation}\label{L74}
   \int_{B_{1+\sqrt{N}}(0)}|v_n|^2\mathrm{d}x> \frac{\delta}{2}.
 \end{equation}
 Since $V(x)$ and $f(x, u)$ are periodic on $x$, we have $\|v_n\|=\|u_n\|$ and
 \begin{equation}\label{L75}
   \Phi(v_n)\rightarrow c_*, \ \ \ \ \|\Phi'(v_n)\|(1+\|v_n\|)\rightarrow 0.
 \end{equation}
 Passing to a subsequence, we have $v_n\rightharpoonup \bar{v}$ in $E$, $v_n\rightarrow \bar{v}$ in $L^{s}_{\mathrm{loc}}(\R^N)$,
 $2\le s<2^*$ and $v_n\rightarrow \bar{v}$ a.e. on $\R^N$. Obviously, \x{L74} implies that $\bar{v}\ne 0$.
 By a standard argument, one has $\Phi'(\bar{v})=0$. This shows that $\bar{v}\in \mathcal{N}^{-}$  and so $\Phi(\bar{v})\ge m$.
 On the other hand, by using \x{L75}, (WN) and Fatou's lemma, we have
 \begin{eqnarray*}
   m & \ge & c_*=\lim_{n\to\infty}\left[\Phi(v_n)-\frac{1}{2}\langle\Phi'(v_n), v_n \rangle\right]
             = \lim_{n\to\infty}\int_{\R^N}\left[\frac{1}{2}f(x, v_n)v_n-F(x, v_n)\right]\mathrm{d}x\\
     & \ge & \int_{\R^N}\lim_{n\to\infty}\left[\frac{1}{2}f(x, v_n)v_n-F(x, v_n)\right]\mathrm{d}x
             =  \int_{\R^N}\left[\frac{1}{2}f(x, \bar{v})\bar{v}-F(x, \bar{v})\right]\mathrm{d}x\\
     &  =  & \Phi(\bar{v})-\frac{1}{2}\langle\Phi'(\bar{v}), \bar{v} \rangle = \Phi(\bar{v}).
 \end{eqnarray*}
 This shows that $\Phi(\bar{v})\le m$ and so $\Phi(\bar{v})=m=\inf_{\mathcal{N}^{-}}\Phi >0$.
 \qed

 \vskip6mm
 {\section{The asymptotically periodic case}}
 \setcounter{equation}{0}

 \vskip4mm
 \par
    In this section, we always assume that $V$ satisfies (V0$'$) and (V1$'$) and define functional $\Phi_0$ as follows:
 \begin{equation}\label{Ph0}
   \Phi_0(u)=\frac{1}{2}\int_{{\R}^N}\left(|\nabla u|^2+V_0(x)u^2\right)\mathrm{d}x-\int_{{\R}^N}F_0(x, u)\mathrm{d}x, \ \ \ \ \forall \ u\in E,
 \end{equation}
 where $F_0(x, t):=\int_{0}^{t}f_0(x, s)\mathrm{d}s$. Then (V0$'$), (F0), (F1) and (F2$'$) imply that $\Phi_0\in C^{1}(E, \R)$ and
 \begin{equation}\label{Phd0}
   \langle \Phi_0'(u), v \rangle  =  \int_{{\R}^N}\left(\nabla u\nabla v+V_0(x)uv\right)\mathrm{d}x-\int_{{\R}^N}f_0(x, u) v\mathrm{d}x,
      \ \ \ \ \forall \ u, v\in E.
 \end{equation}

 \vskip4mm
 \par\noindent
 {\bf Lemma 4.1.}\ \ {\it Suppose that} (V0$'$), (V1$'$), (F0), (F1), (F2$'$), (SQ) {\it and} (WN) {\it are satisfied. Then any sequence
 $\{u_n\}\subset E$ satisfying \x{Ce} is bounded in $E$.}

 \vskip2mm
 \par\noindent
 {\bf Proof.} \ \ To prove the boundedness of $\{u_n\}$, arguing by contradiction, suppose that
 $\|u_n\|_0 \to \infty$. Let $v_n=u_n/\|u_n\|_0$. Then $1=\|v_n\|_0^2$.
 By Sobolev imbedding theorem, there exists a constant $C_4>0$ such that $\|v_n\|_2 \le C_4$.
 Passing to a subsequence, we have $v_n\rightharpoonup \bar{v}$ in $E$. There are two possible cases:
 i). $\bar{v}= 0$ and ii). $\bar{v}\ne 0$.

 \vskip2mm
 \par
   Case i). $\bar{v}=0$, i.e. $v_n\rightharpoonup 0$ in $E$. Then $v_n^{\mathcal{F}+}\rightarrow 0$ and
 $v_n^{\mathcal{F}-}\rightarrow 0$ in $L^{s}_{\mathrm{loc}}(\R^N)$, $2\le s<2^*$
 and $v_n^{\mathcal{F}+}\rightarrow 0$ and $v_n^{\mathcal{F}-}\rightarrow 0$ a.e. on $\R^N$. By (V0$'$), it is easy
 to show that
 \begin{equation}\label{L061}
   \lim_{n\to\infty}\int_{\R^N}V_1(x)(v_n^{\mathcal{F}+})^2\mathrm{d}x=\lim_{n\to\infty}\int_{\R^N}V_1(x)(v_n^{\mathcal{F}-})^2\mathrm{d}x=0.
 \end{equation}

 \par
    If
 $$
   \delta:=\limsup_{n\to\infty}\sup_{y\in \R^N}\int_{B_1(y)}|v_n^{\mathcal{F}+}|^2\mathrm{d}x=0,
 $$
 then by Lions' concentration compactness principle \cite{Lio} or \cite[Lemma 1.21]{Wi}, $v_n^{\mathcal{F}+}\rightarrow 0$ in
 $L^{s}(\R^N)$ for $2<s<2^*$. Fix $R>[2(1+c_*)]^{1/2}$. By virtue of (F0) and (F1), for  $\varepsilon =1/4(RC_4)^2>0$,
 there exists $C_{\varepsilon}>0$ such that \x{F30} holds. Hence, it follows that
 \begin{eqnarray}\label{L063}
   \limsup_{n\to\infty}\int_{\R^N}F(x, Rv_n^{\mathcal{F}+})\mathrm{d}x
     & \le & \limsup_{n\to\infty}\left[\varepsilon R^2\|v_n^{\mathcal{F}+}\|_2^2+C_{\varepsilon}R^p\|v_n^{\mathcal{F}+}\|_p^{p}\right]\nonumber\\
     & \le & \varepsilon(RC_4)^2=\frac{1}{4}.
 \end{eqnarray}
 Let $t_n=R/\|u_n\|_0$. Hence, by virtue of \x{Ce}, \x{L061}, \x{L063} and Corollary 2.5, one can get that
 \begin{eqnarray*}
   c_*+o(1)
   &  =  & \Phi(u_n)\\
   & \ge & \frac{t_n^2}{2}\|u_n\|_0^2-\int_{\R^N}F(x, t_nu_n^{\mathcal{F}+})\mathrm{d}x
             +\frac{1-t_n^2}{2}\langle\Phi'(u_n), u_n \rangle\\
   &     & \ \  +t_n^2\langle\Phi'(u_n), u_n^{\mathcal{F}-} \rangle
             +\frac{t_n^2}{2}\int_{{\R}^N}V_1(x)\left[(u_n^{\mathcal{F}+})^2-(u_n^{\mathcal{F}-})^2\right]\mathrm{d}x\\
   &  =  & \frac{R^2}{2}\|v_n\|_0^2-\int_{\R^N}F(x, Rv_n^{\mathcal{F}+})\mathrm{d}x
             +\left(\frac{1}{2}-\frac{R^2}{2\|u_n\|_0^2}\right)\langle\Phi'(u_n), u_n \rangle\\
   &     & \ \  +\frac{R^2}{\|u_n\|_0^2}\langle\Phi'(u_n), u_n^{\mathcal{F}-} \rangle
            +\frac{R^2}{2}\int_{{\R}^N}V_1(x)\left[(v_n^{\mathcal{F}+})^2-(v_n^{\mathcal{F}-})^2\right]\mathrm{d}x\\
   & \ge & \frac{R^2}{2}-\int_{\R^N}F(x, Rv_n^{\mathcal{F}+})\mathrm{d}x +o(1)\\
   & \ge & \frac{R^2}{2}-\frac{1}{4} +o(1) >  c_*+\frac{3}{4} +o(1).
 \end{eqnarray*}
 This contradiction shows that $\delta>0$.

 \par
    Passing to a subsequence, we may assume the existence of $k_n\in \Z^N$ such that
 $\int_{B_{1+\sqrt{N}}(k_n)}|v_n^{\mathcal{F}+}|^2\mathrm{d}x> \frac{\delta}{2}$. Let $w_n(x)=v_n(x+k_n)$. Since
 $V_0(x)$ is 1-periodic in each of $x_1, x_2, \ldots, x_N$. Then
 \begin{equation}\label{l13}
   \int_{B_{1+\sqrt{N}}(0)}|w_n^{\mathcal{F}+}|^2\mathrm{d}x> \frac{\delta}{2}.
 \end{equation}
 Now we define $\tilde{u}_n(x)=u_n(x+k_n)$, then $\tilde{u}_n/\|u_n\|_0=w_n$ and $\|w_n\|_0=\|v_n\|_0=1$.
 Passing to a subsequence, we have $w_n\rightharpoonup w$ in $E$, $w_n\rightarrow w$ in $L^{s}_{\mathrm{loc}}(\R^N)$,
 $2\le s<2^*$ and $w_n\rightarrow w$ a.e. on $\R^N$. Obviously, \x{l13} implies that $w\ne 0$. Hence, it follows from \x{Ce},
 (SQ) and Fatou's lemma that
 \begin{eqnarray*}
  0 &  =  & \lim_{n\to\infty}\frac{c_*+o(1)}{\|u_n\|_0^2} = \lim_{n\to\infty}\frac{\Phi(u_n)}{\|u_n\|_0^2}\nonumber\\
    &  =  & \lim_{n\to\infty}\left[\frac{1}{2}\left(\|v_n^{\mathcal{F}+}\|_0^2-\|v_n^{\mathcal{F}-}\|_0^2\right)
              +\frac{1}{2}\int_{{\R}^N}V_1(x)\left[(v_n^{\mathcal{F}+})^2-(v_n^{\mathcal{F}-})^2\right]\mathrm{d}x\right.\nonumber\\
    &     & \ \ \left. -\int_{\R^N}\frac{F(x+k_n,\tilde{u}_n)}{\tilde{u}_n^2}w_n^2dx\right]\nonumber\\
    & \le & \frac{1}{2}-\liminf_{n\to\infty}\int_{\R^N}\frac{F(x+k_n, \tilde{u}_n)}{\tilde{u}_n^2}w_n^2dx
             \le \frac{1}{2}-\int_{\R^N}\liminf_{n\to\infty}\frac{F(x+k_n, \tilde{u}_n)}{\tilde{u}_n^2}w_n^2dx\nonumber\\
    &  =  & -\infty,
 \end{eqnarray*}
 which is a contradiction.

 \vskip2mm
 \par
   Case ii). $\bar{v}\ne 0$. In this case, we can also deduce a contradiction by a standard argument.

 \vskip2mm
 \par
 Cases i) and ii) show that $\{u_n\}$ is bounded in $E$.
 \qed

 \vskip4mm
 \par\noindent
 {\bf Proof of Theorem 1.3.} \ \ Applying Lemmas 2.10 and 4.1, we deduce that there exists a bounded sequence
 $\{u_n\}\subset E$ satisfying \x{Ce}. Passing to a subsequence, we have $u_n\rightharpoonup \bar{u}$ in $E$.
 Next, we prove $\bar{u}\ne 0$.

 \par
   Arguing by contradiction, suppose that $\bar{u}=0$, i.e. $u_n\rightharpoonup 0$ in $E$, and so
 $u_n\rightarrow 0$ in $L^{s}_{\mathrm{loc}}(\R^N)$, $2\le s<2^*$ and $u_n\rightarrow 0$ a.e. on $\R^N$. By (V0$'$) and (F2$'$),
 it is easy to show that
 \begin{equation}\label{T32}
   \lim_{n\to\infty}\int_{\R^N}V_1(x)u_n^2dx=0, \ \ \ \
   \lim_{n\to\infty}\int_{\R^N}V_1(x)u_nvdx=0, \ \ \ \ \forall \ v\in E
 \end{equation}
 and
 \begin{equation}\label{T34}
   \lim_{n\to\infty}\int_{\R^N}F_1(x, u_n)dx=0, \ \ \ \ \lim_{n\to\infty}\int_{\R^N}f_1(x, u_n)vdx=0, \ \ \ \ \forall \ v\in E.
 \end{equation}
 Note that
 \begin{equation}\label{T35}
   \Phi_0(u)=\Phi(u)-\frac{1}{2}\int_{\R^N}V_1(x)u^2dx+\int_{\R^N}F_1(x, u)dx, \ \ \ \ \forall \ u\in E
 \end{equation}
 and
 \begin{equation}\label{T36}
   \langle\Phi_0'(u), v\rangle=\langle\Phi'(u), v\rangle-\int_{\R^N}V_1(x)uvdx+\int_{\R^N}f_1(x, u)vdx, \ \ \ \ \forall \ u, v\in E.
 \end{equation}
 From \x{Ce}, \x{T32}-\x{T36}, one can get that
 \begin{equation}\label{T37}
   \Phi_0(u_n)\rightarrow c_*, \ \ \ \ \|\Phi_0'(u_n)\|(1+\|u_n\|_0)\rightarrow 0.
 \end{equation}

 \par
    Analogous to the proof of Theorem 1.2, we may assume the existence of $k_n\in \Z^N$ such that \\
 $\int_{B_{1+\sqrt{N}}(k_n)}|u_n|^2\mathrm{d}x > \frac{\delta}{2}$ for some $\delta>0$. Let $v_n(x)=u_n(x+k_n)$. Then
 $\|v_n\|_0=\|u_n\|_0$ and
 \begin{equation}\label{T38}
   \int_{B_{1+\sqrt{N}}(0)}|v_n|^2dx> \frac{\delta}{2}.
 \end{equation}
 Passing to a subsequence, we have $v_n\rightharpoonup \bar{v}$ in $E$, $v_n\rightarrow \bar{v}$ in $L^{s}_{\mathrm{loc}}(\R^N)$,
 $2\le s<2^*$ and $v_n\rightarrow \bar{v}$ a.e. on $\R^N$. Obviously, \x{T38} implies that $\bar{v}\ne 0$.
 Since $V_0(x)$ and $f_0(x, u)$ are periodic in $x$, then by \x{T37}, we have
 \begin{equation}\label{T39}
   \Phi_0(v_n)\rightarrow c_*, \ \ \ \ \|\Phi_0'(v_n)\|(1+\|v_n\|_0)\rightarrow 0.
 \end{equation}
 In the same way as the last part of the proof of Theorem 1.2, we can prove that $\Phi_0'(\bar{v})=0$
 and $\Phi_0(\bar{v})\le c_*$.

 \par
   It follows from $\Phi_0'(\bar{v})=0$ and \x{Phd0} that $\bar{v}^{\mathcal{F}+}\ne 0$. By Lemma 2.11, there
 exist $t_0=t(\bar{v})>0$ and $w_0=w(\bar{v})\in E^{\mathcal{F}-}$ such that $t_0\bar{v}+w_0\in \mathcal{N}^{0}$,
 and so $\Phi(t_0\bar{v}+w_0)\ge m$.  By virtue of (F2$'$), $f_0(x, t)/|t|$ is non-decreasing on
 $t\in (-\infty, 0)\cup (0, \infty)$, similar to \x{L32}, we have
 \begin{equation}\label{T39}
   \frac{1-t_0^2}{2}f_0(x, \bar{v})\bar{v}-t_0f_0(x, \bar{v})w_0-\int_{t_0\bar{v}+w_0}^{v}f_0(x, s)\mathrm{d}s\ge 0.
 \end{equation}
 Hence, from \x{Ph0}, \x{Phd0}, \x{T39} and the fact that $-V_1(x)t^2+F_1(x, t)>0$ for
 $|x|<1+\sqrt{N}$ and $t\ne 0$, we have
 \begin{eqnarray*}
  m & \ge & c_*\ge \Phi_0(\bar{v})\\
    &  =  & \Phi_0(t_0\bar{v}+w_0)+\frac{1}{2}\|w_0\|_0^2+\frac{1-t_0^2}{2}\langle\Phi_0'(\bar{v}), \bar{v} \rangle
              -t_0\langle\Phi_0'(\bar{v}), w_0 \rangle\\
    &     &   +\int_{\R^N}\left[\frac{1-t_0^2}{2}f_0(x, \bar{v})\bar{v}-t_0f_0(x, \bar{v})w_0
              -\int_{t_0\bar{v}+w_0}^{\bar{v}}f_0(x, s)\mathrm{d}s\right]\mathrm{d}x\\
    & \ge & \Phi_0(t_0\bar{v}+w_0)+\frac{1}{2}\|w_0\|_0^2\\
    &  =  & \frac{1}{2}\|w_0\|_0^2+\Phi(t_0\bar{v}+w_0)-\frac{1}{2}\int_{\R^N}V_1(x)(t_0\bar{v}+w_0)^2\mathrm{d}x
             +\int_{\R^N}F_1(x, t_0\bar{v}+w_0)\mathrm{d}x\\
    &  >  & \Phi(t_0\bar{v}+w_0)\ge m,
 \end{eqnarray*}
 since $\bar{v}(x)\not\equiv 0$ for $x\in B_{1+\sqrt{N}}(0)$. This contradiction implies that $\bar{u}\ne 0$.
 In the same way as the last part of the proof of Theorem 1.2, we can certify that $\Phi'(\bar{u})=0$ and
 $\Phi(\bar{u})=m=\inf_{\mathcal{N}^{0}}\Phi$. This shows that $\bar{u}\in E$ is a solution for problem \x{SP} with
 $\Phi(\bar{u})=\inf_{\mathcal{N}^{0}}\Phi>0$.
 \qed

 {}
\end{document}